\documentclass[12pt]{article}

\usepackage{amsfonts}

\newtheorem{theorem}{Theorem}[section]
\newtheorem{lemma}[theorem]{Lemma}
\newtheorem{corollary}[theorem]{Corollary}
\newtheorem{definition}[theorem]{Definition}
\newtheorem{proposition}[theorem]{Proposition}
\newtheorem{example}[theorem]{Example}
\newtheorem{remark}[theorem]{Remark}

\def\bit{\begin{itemize}}
\def\eit{\end{itemize}}

\def\bc{\begin{center}}
\def\ec{\end{center}}

\def\bthm{\begin{theorem}}
\def\ethm{\end{theorem}}

\def\bcor{\begin{corollary}}
\def\ecor{\end{corollary}}

\def\bprop{\begin{proposition}}
\def\eprop{\end{proposition}}

\def\blem{\begin{lemma}}
\def\elem{\end{lemma}}

\def\bex{\begin{example}}
\def\eex{\end{example}}

\def\brem{\begin{remark}}
\def\erem{\end{remark}}

\def\prf{\noindent{\bf Proof. }}

\def\bdes{\begin{description}}
\def\edes{\end{description}}

\def\ita{\item[(a)]}
\def\itb{\item[(b)]}
\def\itc{\item[(c)]}

\def\iti{\item[(i)]}
\def\itii{\item[(ii)]}

\def\beq{\begin{equation}}
\def\eeq{\end{equation}}

\def\ben{\begin{enumerate}}
\def\een{\end{enumerate}}

\def\beqar{\begin{eqnarray}}
\def\eeqar{\end{eqnarray}}

\def\beqarr{\begin{eqnarray*}}
\def\eeqarr{\end{eqnarray*}}

\def \non{{\nonumber}}

\def\RR{{\mathbb R}}  

\def\cE{\mathcal{E}}
\def\cF{\mathcal{F}}

\def\cM{\mathcal{M}}
\def\cP{\mathcal{P}}

\def\qed{\vbox{\hrule\hbox{\vrule height 1.5 ex\kern 1 ex\vrule}\hrule}}

\def\P{{\mathsf P}} 

\def\E{{\mathsf E}} 

\def\NN{{\mathbb N}}       

\def\eps{\varepsilon}

\def\p{\varphi}

\def\la{\langle}
\def\ra{\rangle}

\def\part{\partial}

\begin{document}

\title{Sticky flows on the circle.}
\author{Yves Le Jan and Olivier Raimond}
\maketitle


\section*{Introduction}
The purpose of this note is to give an example of stochastic flows of
kernels as defined in \cite{ljr}, which
naturally interpolates between the Arratia coalescing flow associated
with systems of coalescing independent Brownian particles on the
circle and the deterministic diffusion flow (actually, the results are
given in the slightly more general framework of symmetric Levy
processes for which points are not polar). The construction is
performed using Dirichlet form theory and the extension of De
Finetti's theorem given in \cite{ljr}. The sticky flows of kernels are
associated with systems of sticky independent Levy particles on the
circle, for some fixed parameter of stickyness. Some elementary
asymptotic properties of the flow are also given.

\section{Compatible family of Dirichlet forms.}
Let $(\cE_n)_{n\geq 1}$ be a family of Dirichlet forms\footnote{we
refer the reader not familiar with Dirichlet forms and symmetric
Markov processes to \cite{fuku}}, respectively
defined on $L^2(M^n,m_n)$, where $M$ is a metric space and
$(m_n)_{n\geq 1}$ is a family of probability measures on
$M^n$. We will denote by $D_n$ the domain of the Dirichlet form
$\cE_n$. For all $n\geq 1$, $S_n$ denotes the group of permutations of
$\{1,\dots,n\}$.

Let $(\P^{(n)}_t)_{n\geq 1}$ be the family of Markovian semigroups
associated with this family of Dirichlet forms.
\begin{definition} We will say that the family of Dirichlet forms
 $(\cE_n)_{n\geq 1}$ is compatible if the family of Markovian
 semigroups $(\P^{(n)}_t)_{n\geq 1}$ is compatible,
that is if the following assertions are satisfied
\bdes
\iti for all $f\in L^2(m_n)$ and $\sigma\in S_n$,
$\P^{(n)}_t f_\sigma = (\P^{(n)}_t f)_\sigma$,
where $f_\sigma(x_1,\dots,x_n)=f(x_{\sigma_1},\dots,x_{\sigma_n})$;
\itii for all $f\in L^2(m_n)$, $f\otimes 1\in L^2(m_{n+1})$ and
$\P^{(n+1)}_t(f\otimes 1)  = (\P^{(n)}_tf)\otimes 1$.
\edes
\end{definition}

\begin{definition} We will say that the family of probability measures
 $(m_n)_{n\geq 1}$ is consistent and exchangeable if for
  all $n\geq 1$, $m_n$ is the law of $(X_1,\dots,X_n)$, where
  $(X_i)_{i\geq 1}$ is an exchangeable sequence of $M$-valued random
  variables. By Kolmogorov's theorem, this holds if and only if
\bdes \iti for all $n\geq 1$, $\sigma\in S_n$ and $f\in L^1(m_n)$,
$\int f_\sigma dm_n=\int f dm_n$.
\itii for all $n\geq 1$ and $f\in L^1(m_n)$,
$\int (f\otimes 1) dm_{n+1}=\int f dm_n$.
\edes
\end{definition}

We now assume that $(m_n)_{n\geq 1}$ is consistent and exchangeable (note that
this holds if the family of Dirichlet forms is compatible and if they
are irreducible). Then we can define for all $n\geq 1$ the conditional
expectations $\pi_n:L^2(m_{n+1})\to L^2(m_n)$ such that for all $f\in
L^2(m_{n+1})$, $\pi_n(f)\otimes 1$ is the orthogonal projection in
$L^2(m_{n+1})$ of $f$ onto $\{g\otimes 1,~g\in L^2(m_{n})\}$.

\bprop\label{th:comp}
The family of Dirichlet forms $(\cE_n)_{n\geq 1}$ is compatible
if and only if
\bdes \iti For all $n\geq 1$, $\sigma\in S_n$ and $(f,g)\in
D_n^2$, we have $(f_\sigma,g_\sigma)\in D_n^2$ and
$\cE_n(f_\sigma,g_\sigma) = \cE_n(f,g)$.
\itii For all $n\geq 1$, $f\in D_{n+1}$ and $g\in D_n$, we have
$g\otimes 1\in D_{n+1}$, $\pi_n(f)\in D_n$ and
\begin{eqnarray}
&&\cE_{n+1}(f,g\otimes 1) \quad = \quad \cE_n(\pi_n(f),g),\\
&& \cE_{n+1}(f) \quad = \quad \cE_n(\pi_n(f)) + \cE_{n+1}(f-\pi_n(f)\otimes 1).\label{decomp}
\end{eqnarray}
\edes
\eprop
This proposition is an immediate corollary of the following theorem.

\bthm Let $\cE_1$ and $\cE$ be two Dirichlet forms respectively
defined on $L^2(E_1,\cF_1,m_1)$ and on 
$L^2(E_1\times E_2,\cF_1\otimes \cF_2,m)$, with domain $D_1$ and $D$,
where $m_1$ and $m$ are probability measures on $E_1$ and on
$E_1\times E_2$ satisfying $\int (g\otimes 1)dm = \int gdm_1$ for all
$f\in L^1(m_1)$. Let $\P^1_t$ and $\P_t$ be the associated Markovian
semigroups. Then, {\bf (i)} and {\bf (ii)} are equivalent, where
\bdes \iti For all $g\in L^2(m_1)$, 
$\P_t(g\otimes 1)=(\P^1_tg)\otimes 1$.
\itii For all $g\in D_1$ and $f\in D$, we have $\pi(f)\in D_1$,
$g\otimes 1\in D$ and 
\begin{eqnarray} \label{proj} && \cE(f,g\otimes 1) \quad = \quad \cE_1(\pi(f),g),\\
&& \cE(f) \quad = \quad \cE_1(\pi(f)) + \cE(f-\pi(f)\otimes 1). \end{eqnarray}
where $\pi(f)\otimes 1$ is the orthogonal projection of $f$ in
$L^2(m)$ onto $\{g\otimes 1,~g\in L^2(m_1)\}$. 
\edes \ethm
\prf Suppose first {\bf (i)}.

For all $f\in L^2(m)$, $\cE(f)=\lim_{t\to  0}
\left(\|f\|^2_{L^2(m)}-\|\P_{t/2}f\|^2_{L^2(m)}\right)$ and $f\in D$
if and only if $\cE(f)<\infty$. Using this relation and (i), we show
that for all $g\in L^2(m_1)$, $\cE(g\otimes 1)=\cE_1(g)$. Thus $g\in
D_1$ if and only if $g\otimes 1\in D$.

Let $f\in L^2(m)$ and $u=f-\pi(f)\otimes 1$. Then $u$ is orthogonal in
$L^2(m)$ to $V=\{g\otimes 1,~g\in L^2(m_1)\}$. Then for all positive
$t$, $\P_tu$ is orthogonal to $V$ since 
$$\la \P_t u,g\otimes
1\ra_{L^2(m)} =  \la u,\P_t(g\otimes 1)\ra_{L^2(m)} = \la
u,(\P^1_tg)\otimes 1\ra_{L^2(m)} = 0.$$
This implies that $\|\P_{t/2}f\|^2_{L^2(m)}=\|\P_{t/2}u\|^2_{L^2(m)} + 
\|\P_{t/2}^1\pi(f)\|^2_{L^2(m_1)}$. Moreover,
$\|f\|^2_{L^2(m)}=\|u\|^2_{L^2(m)}+\|\pi(f)\|^2_{L^2(m_1)}$ and we
prove $\cE(f)=\cE(u)+\cE_1(\pi(f))$. It follows that
$\cE_1(\pi(f))\leq\cE(f)$, which implies that $\pi(f)\in D_1$ when
$f\in D$.

To prove {\bf(ii)}, it remains to prove (\ref{proj}). This follows
from the computation
\beqarr
\cE(f,g\otimes 1) &=& \lim_{t\to 0} 
\frac{1}{t}(\la f,g\otimes 1\ra_{L^2(m)} - 
\la f,\P_t(g\otimes 1)\ra_{L^2(m)})\\
&=& \lim_{t\to 0} 
\frac{1}{t}(\la \pi(f),g\ra_{L^2(m_1)} - 
\la \pi(f),\P^1_t g\ra_{L^2(m_1)})\\
&=& \cE_1(\pi(f),g). \eeqarr

We now assume {\bf (ii)} is satisfied. Then {\bf (i)} is satisfied if
for all $\alpha>0$ and $g\in L^2(m_1)$, 
$G_\alpha(g\otimes 1)=(G^1_\alpha g)\otimes 1$, where
$G_\alpha$ and $G^1_\alpha$ are respectively the resolvents of $\P_t$
and of $\P^1_t$. Let $\cE_\alpha$ and $\cE^1_\alpha$ be respectively
the forms $\cE+\alpha\la\cdot,\cdot\ra_{L^2(m)}$  and
$\cE_1+\alpha\la\cdot,\cdot\ra_{L^2(m_1)}$. Let $g\in L^2(m_1)$ and
$\alpha>0$, then $G_\alpha(g\otimes 1)$ is the unique element of
$L^2(m)$ such that for all $f\in L^2(m)$ we have
$\cE_\alpha(f,G_\alpha(g\otimes 1))=\la f,g\otimes
1\ra_{L^2(m)}$. This implies that 
$$\cE_\alpha(f,G_\alpha(g\otimes
1))=\la \pi(f),g\ra_{L^2(m_1)}.$$

Using equation (\ref{proj}) and the definition of $G^1_\alpha g$, we
also have 
$$\cE_\alpha(f,(G^1_\alpha g)\otimes 1) = 
\cE^1_\alpha(\pi(f),G^1_\alpha g) = 
\la \pi(f),g\ra_{L^2(m_1)}.$$
This proves that
$G_\alpha (g\otimes 1)=(G^1_\alpha g)\otimes 1$. \qed

\section{Compatible families of probability measures and exchangeable random partitions.}
For generalities on exchangeable random partitions, we refer to
Pitman's St Flour course \cite{pitman}. 
For all $n\geq 1$, we let $\cP_n$ denote the set of all partitions of
$[n]=\{1,\dots,n\}$. The number of elements of a partition $\pi$ is
denoted $|\pi|$. A random partition $\Pi_n$ of $[n]$ is called
{\em exchangeable} if its distribution is invariant under the obvious
action of $S_n$ on $\cP_n$. Equivalently, for each partition
$\{A_1,\dots,A_k\}$ of $[n]$,
$$\P[\Pi_n=\{A_1,\dots,A_k\}]=p(|A_1|,\dots,|A_k|)$$
for some symmetric function $p$, called the {\em exchangeable partition probability function}
(EPPF) of $\Pi_n$.

A sequence of exchangeable random partitions $(\Pi_n)_{n\geq 1}$ is
called {\em consistent in distribution} if for all $1\leq m\leq n$, the
restriction, denoted by $\Pi_{m,n}$, of $\Pi_n$ to $[m]$ has the same
distribution as $\Pi_m$. The associated sequence of EPPFs are also
called {\em{consistent}}. Note that being given a consistent sequence of
EPPFs $(p_n)$, it is possible to construct a sequence of exchangeable random
partitions, called an {\em{infinite exchangeable random partition}}
 $\Pi_\infty=(\Pi_n)_{n\geq 1}$ such that the EPPF of
$\Pi_n$ is $p_n$ and with $\Pi_{m,n}=\Pi_m$ for all $1\leq m\leq n$.

Using Kingman's representation theorem, given an infinite exchangeable
random partition $\Pi_\infty$, it is possible to construct a random
sequence in $[0,1]$, $(P_i)_{i\geq 1}$, with $\sum_{i\geq 1}P_i\leq
1$ and the law of $\Pi_\infty$ given $(P_i)_{i\geq 1}$ is the same as
if $\Pi_\infty$ were given by random sampling from a random distribution
with ranked atoms $(P_i)_{i\geq 1}$. The random partition $\Pi_\infty$
is called {\em{proper}} when $\sum_i P_i=1$.

For all partition $\pi=\{A_1,\dots,A_k\}$ of $[n]$, let $E_\pi$ denote the
set of all $x\in M^n$ such that for all $1\leq l\leq k$ and  $(i,j)\in
A_l^2$, we have $x_i=x_j$. Then $E_\pi$ is isomorphic to $M^k$. Let
$\lambda_\pi$ be the probability measure on $E_\pi$ given by 
$(\p_\pi)_*(\lambda^{\otimes k})$, where $\lambda$ is a probability
measure on $M$ and $\p_\pi:M^k\to E_\pi$ with
$(\p_\pi(y_1,\dots,y_k))_i=y_j$ for all $i,j$ such that $i\in A_j$.

Let us be given a proper infinite random partition $\Pi_\infty$.
Let $m_n=\sum_{\pi\in\cP_n} p_{\pi}\lambda_\pi$, where
$p_n=(p_\pi)_{\pi\in\cP_n}$ is the distribution of $\Pi_n$, the
restriction of $\Pi_\infty$ to $[n]$. Then $(m_n)_{n\geq 1}$ is
consistent and exchangeable.

Kingman's representation theorem then implies that
$m_n=\E[\mu^{\otimes n}]$, where $\mu$ is a random measure on
$M$. This random measure can be described with $(X_i)_{i\geq 1}$ a
sequence of independent $M$-valued random variables of law $\lambda$
and an independent sequence, $(P_i)_{i\geq 1}$, of $[0,1]$-valued
random variables with $\sum_{i\geq 1}P_i= 1$, and $\mu$ is defined by
the relation  $\mu=\sum_{i\geq 1}P_i\delta_{X_i}$. 

In this paper we will be interested in the case where $(P_i)_{i\geq
1}$ is distributed like a Dirichlet process of parameter
$\frac{1-\tau}{\tau}$, where $\tau\in [0,1]$ is a fixed parameter.
More precisely, the sequence $(P_i)_{i\geq 1}$ is distributed as the
jumps of a process
$\left(\frac{\Gamma_{u\theta}}{\Gamma_\theta}\right)_{0\leq u\leq 1}$,
where $(\Gamma_s)_{s\geq 0}$ is a standard $\Gamma$-process, i.e. the
subordinator whose marginal laws are given by gamma distribution of
parameter $s$. In this case the family $(m_n)_{n\geq 1}$ can be
constructed by the relation $m_1=\lambda$ and $m_{n+1}(d\bar{x}_n,dx_{n+1}) = 
m_{n}(d\bar{x}_n)\pi_n(\bar{x}_n,dx_{n+1})$
where $\bar{x}_n=(x_1,\dots,x_n)$ and 
$$\pi_n(\bar{x}_n,dx_{n+1}) =
\frac{(1-\tau )\lambda(dx_{n+1})
+\tau\sum_{i=1}^n\delta_{x_i}(dx_{n+1})}{(1-\tau)+n\tau}.$$
Then the EPPF $p$ satisfies
\beq \label{eppf}\left\{\begin{array}{lll} p(n_1,\dots,n_k,1) &=&
\left(\frac{1-\tau}{(1-\tau)+n\tau}\right)p(n_1,\dots,n_k),\\
p(n_1,\dots,n_k+1) &=&
\left(\frac{n_k\tau}{(1-\tau)+n\tau}\right)p(n_1,\dots,n_k).
\end{array}\right.\eeq

\section{Compatible families of Dirichlet forms on $S^1$ and sticky flows.}
We now let $M$ denote the unit circle $S^1$ and $\lambda$ denote the
Lebesgue measure on $S^1$.
For a fixed parameter $\tau\in [0,1]$, assume we are given the
compatible family of probability measures $(m_n)_{n\geq 1}$ defined in
the previous section. Let $\P_t$ be the Markovian semigroup
of a symmetric Levy process of exponent $\psi$ on $S^1$, for which
points are not polar, i.e. such that
$\int_{S^1}\left(\frac{1}{\alpha+\psi(x)}\right) dx<\infty$ for
$\alpha>0$. We denote by $\cE$ the associated Dirichlet form. This
Dirichlet form is defined on $L^2(\lambda)$. 

For all $n\geq 1$, we define the Dirichlet form on
$C^1((S^1)^n)\subset L^2(m_n)$ by the formula
\beq\cE_n=\sum_{\pi\in\cP_n}p_\pi\cE_\pi, \label{defEn}\eeq
where $\cE_\pi=\p_\pi(\cE^{\odot k})$ (with $k=|\pi|$). More
precisely, for $f\in C^1(E_\pi)$, there exists $g\in C^1((S^1)^k)$
with $f=g\circ \p_\pi$, then $\cE_\pi(f)=\cE^{\odot k}(g)$. Here
$\cE^{\odot k}$ denotes the Dirichlet form associated with $k$
independent Levy processes, i.e. with the Markovian semigroup
$\P_t^{\otimes k}$. 

\bprop
This family of Dirichlet forms verifies the following recurrence
property~:
\beq\label{recurs} \cE_{n+1}(g)
=\frac{1}{(1-\tau)+n\tau}\left((1-\tau) \cE_n \odot \cE(g)
+\tau \sum_{i=1}^n \cE_n(g^i)\right),\eeq
where $g^i(x_1,\dots,x_n)$ denotes $g(x_1,\dots,x_n,x_i)$ and $\cE_n
\odot \cE$ is the Dirichlet form defined on 
$C^1((S^1)^n)\subset L^2(m_n\otimes\lambda)$ associated with the
Markovian semigroup $\P^{(n)}_t\otimes\P_t$. \eprop 

\prf We have
$$ \cE_{n+1}(g) = \sum_{\pi\in\cP_{n}}p_{\pi\cup\{n+1\}} \cE_{\pi\cup\{n+1\}}(g)
+ \sum_{\pi\in\cP_{n}}\sum_{i=1}^n p_{\pi_i} \cE_{\pi_i}(g), $$
where $\pi_i$ is the partition of $[n+1]$ obtained by adding $n+1$ to the
set $\pi$ containing $i$. We conclude using the definition of the EPPF
$p$ given by (\ref{eppf}), the fact that $\cE_{\pi\cup\{n+1\}}(g) =
\cE_{\pi}\odot \cE(g)$ and that $\cE_{\pi_i}(g)=\cE_\pi(g^i)$. \qed

\medskip
The Dirichlet form $\cE_n$ being defined by the superposition of closable
Dirichlet forms, $\cE_n$ is closable in $L^2(m_n)$ (see proposition
3.1.1 p.214 in \cite{bouleau}). Moreover, these Dirichlet forms are
regular by construction.

\brem Let $A_n$ denote the generator of $\cE_n$ and let $A^{(n)}$
denote the generator of  $\P_t^{\otimes n}$. (When $\P_t$ is the heat
semigroup, $A^{(n)}$ is the Laplacian on $H_2((S^1)^n)$). Then for all
$\alpha>0$, $\{(-A_n+\alpha)f,~f\in C^\infty((S^1)^n)\}$ is dense in
$L^2(m_n)$ and we have 
\beq A_nf~m_n = \sum_{\pi\in\cP_n}p_\pi 
(A^{(|\pi|)}(f\circ \p_\pi))\circ \p_\pi^{-1}~\lambda_\pi. \eeq
\erem

\bthm \bdes 
\iti The family of Dirichlet forms $(\cE_n)_{n\geq 1}$ defined above is
compatible. 
\itii The family of Markovian semigroups $(\P^{(n)}_t)_{n\geq 1}$
associated with this family of Dirichlet forms are strong Feller semigroups. 
\edes \ethm

\prf {\bf(i)}: For all $n\geq 1$, $C^\infty((S^1)^n)$ is dense in $D_n$,
the domain of $\cE_n$, with respect to the scalar product
$\cE_n+\la\cdot,\cdot\ra_{L^2(m_n)}$. From the definition of
$(m_n)_{n\geq 1}$, the projection operator $\pi_n$ maps
$C^{\infty}((S^1)^{n+1})$ onto $C^{\infty}((S^1)^n)$. Since
$\cE_n(\pi_n(f))\leq \cE_{n+1}(f)$ (see (\ref{decomp})) and
$\|\pi_n(f)\|_{L^2(m_{n+1})}\leq\|f\|_{L^2(m_n)}$, $\pi_n$ maps
$D_{n+1}$ onto $D_n$.

Let now $g\in D_n$, then using (\ref{recurs}) it is easy to check that 
$\cE_{n+1}(g\otimes 1) = \cE_n(g)$, which implies $g\otimes 1 \in D_{n+1}$.

Let us remark that 
$$\cE^{\odot (k+1)}(f,g\otimes 1) = 
\int \cE^{\odot k}(f_x,g) \lambda(dx),$$
 where $f_x(x_1,\dots,x_k)=f(x_1,\dots,x_k,x)$. This holds since
$$\cE^{\odot (k+1)}(f,g\otimes 1) = \lim_{t\to 0}\frac{1}{t}\left(\la
f,g\otimes 1\ra_{L^2(\lambda^{\otimes (k+1)})} - \la
f,(\P_t^{\otimes k}g)\otimes 1\ra_{L^2(\lambda^{\otimes (k+1)})}\right).$$
This implies, using (\ref{defEn}), that $\cE_n\odot \cE(f,g\otimes 1) = \int \cE_n(f_x,g)
\lambda(dx)$. Using this relation and (\ref{recurs}), we show that for
all $n\geq 1$, all $f\in D_{n+1}$ and all $g\in D_n$, we have
$\cE_{n+1}(f,g\otimes 1) =\cE_n(\pi_n(f),g)$. 
Thus, we conclude applying theorem \ref{th:comp}
that the family of Dirichlet forms is compatible.

\medskip
{\bf (ii)}: 
We now define an everywhere defined version $\tilde{\P}^{(n)}_t$ of
$\P^{(n)}_t$.

Fix $n\geq 1$. Let $\hat{X}$ be the stationary Markov process
associated with $\P^{(n+1)}_t$. The stationary law is $m_{n+1}$. Let
$X_t$ and $Y_t$ be the processes defined by
$$\begin{array}{lllll} X^i_t &=& \hat{X}^i_t & \mbox{for} & i\leq n,\\
Y^i_t &=& \hat{X}^i_t & \mbox{for} & i\leq n-1,\\
Y^n_t &=& \hat{X}^{n+1}_t & \mbox{for} & t\leq T,\\
Y^n_t &=& \hat{X}^{n}_t & \mbox{for} &  t> T,
\end{array}$$
where $T = \inf\{s,~\hat{X}^n_s=\hat{X}^{n+1}_s\}$.
Note that for $t\geq T$, $X_t=Y_t$. Moreover, the compatibility of the
family $(\P^{(n)}_t)_{n\geq 1}$ and the strong Markov property implies
that $(X_t)$ and $(Y_t)$ are both stationary Markov processes
associated with $\P^{(n)}_t$.

For any bounded function $f$ on $(S^1)^n$ and $t>0$,
and any bounded function $h$ on $(S^1)^{n+1}$, we have (in the
following, $x=(\hat{x}_1,\dots,\hat{x}_n)$ and
$y=(\hat{x}_1,\dots,\hat{x}_{n-1},\hat{x}_{n+1})$)
 using the compatibility of the family $(\P^{(n)}_t)_{n\geq 1}$ that 
\beqar
\hskip-30pt\left|\int_{(S^1)^{n+1}} (\P^{(n)}_tf(x)-\P^{(n)}_tf(y))h(\hat{x})
~ m_{n+1}(d\hat{x})\right|\hskip-110pt && \hskip90pt =\non\\ 
&=& |\E[(f(X_t)-f(Y_t))h(\hat{X}_0)]|\non\\
&\leq& 2\|f\|_\infty \E[1_{t\geq T}|h(\hat{X}_0)|]\non\\
&\leq& 2\|f\|_\infty \int |h(\hat{x})|
\eps_t(d(\hat{x}_n,\hat{x}_{n+1})) ~m_{n+1}(d\hat{x}),\label{estim}
\eeqar
where $\eps_t(r)$ is the probability that two independent Levy
processes of exponent $\psi$, respectively started at $x$ and
$y$ with $r=d(x,y)$, have not met before time $t$ .
Therefore, since points are not polar for the Levy process, for all
$t>0$, $\lim_{r\to 0}\eps_t(r)=0$. 

Since the estimate (\ref{estim}) holds for all bounded function $h$,
we have
\beq |\P^{(n)}_tf(x)-\P^{(n)}_tf(y)| \leq
\eps_t(d(\hat{x}_n,\hat{x}_{n+1})) \qquad
m_{n+1}(d\hat{x})-a.e.\eeq
Thus, this implies that
\beq |\P^{(n)}_tf(x)-\P^{(n)}_tf(y)| \leq
\sum_{i=1}^n\eps_t(d(x_i,x_j)) \qquad
m_{2n}(dx,dy)-a.e.\eeq
Therefore there exists a continuous function $\tilde{\P}^{(n)}_tf$
such that $\tilde{\P}^{(n)}_tf=\P^{(n)}_tf$ $m_n(dx)$-a.e.
Thus, this defines an everywhere defined version of $\P^{(n)}_t$,
which is a strong Feller semigroup. \qed

\medskip Let us now apply theorem 1.1.4 in \cite{ljr}. This
theorem (a generalization of De Finetti's theorem) states that from a
compatible family of Feller semigroups it is possible to construct a
stochastic flow of kernels $(K_{s,t},~s\leq t)$ such that
$\E[K_{0,t}^{\otimes n}]=\P^{(n)}_t$. Thus to the family of Dirichlet
forms $(\cE_n)_{n\geq 1}$ is associated a stochastic flow of kernels
which we will call the sticky flow of parameter $\tau$ and exponent $\psi$.

\bprop Let $(K_{s,t},~s\leq t)$ be a sticky flow of parameter $\tau$
and exponent $\psi$. Let $\mu$ denote an independent random probability
measure $\sum_{i\geq 1} P_i\delta_{X_i}$, where $(P_i)$ is a Dirichlet
process of parameter $\frac{1-\tau}{\tau}$ and $(X_i)$ is a sequence
of independent random variables of law $\lambda$ and independent of
$(P_i)$. Then,
\bdes \ita For all $t\in\RR$, $\lambda K_{-T,t}$ converges weakly as
$T\to\infty$ towards a probability measure on $S^1$, denoted
$\mu_t$. Moreover, $\mu_t$ has the same law as $\mu$.
\itb The sticky flow induces a Feller process $\nu_t$ on the space
$\cM_1^+(S^1)$ of probability measures on $S^1$ by the relation
$\nu_t=\nu_0K_{0,t}$. The stationary distribution of this Feller
process is given by the law of $\mu$.
\itc For all $x\in S^1$ and all $s\leq t$, a.s. $K_{s,t}(x)$ is an atomic
measure.
\edes
\eprop
\prf {\bf (a) and (b)}~: The fact that $\lambda K_{-T,t}$ converges
weakly follows from the fact that 
$\lambda K_{-T,t} f=\int K_{-T,t}f(x)\lambda(dx)$ is a martingale in
$T$, which therefore converges a.s. The Feller property of the
semigroup associated with $\nu_t$ can be easily proved from the Feller
properties of the semigroups $\P^{(n)}_t$ using the dense algebra of
polynomial functions on $\cM_1^+(S^1)$ of the form
$$\hat{g}(\nu)=\int g(x_1,\dots,x_n)~\nu(dx_1)\cdots\nu(dx_n), \qquad
g\in C((S^1)^n),\quad n\in\NN.$$

It is clear that the law of $\mu_t$ is a stationary distribution
for this Feller process on $\cM(S^1)$. Thus it
remains to show that the law of $\mu$ is the unique stationary
distribution. 

Note that for all $n\geq 1$, 
$\E[(\nu_t)^{\otimes n}]=\nu_0^{\otimes n}\P^{(n)}_t$. 
Since $\P^{(n)}_t$ is irreducible, $\E[(\nu_t)^{\otimes n}]$
converges towards $m_n=\E[\mu^{\otimes n}]$. This implies that if
the law of $\nu$ is a stationary distribution, then for all $n\geq 1$
we have $\E[\nu^{\otimes n}]=\E[\mu^{\otimes n}]$. This implies, since
$\cM_1^+(S^1)$ is compact, that $\nu$ and $\mu$ have the same law. 

{\bf (c)} Since $\mu$ is atomic and stationary, $\mu K_{0,t}$ is
atomic and distributed like $\mu$. Since $\mu K_{0,t} = \int
K_{0,t}(x) \mu(dx)$, we have $\mu(dx)$-a.e. $K_{0,t}(x)$ is
atomic a.s. Thus, since $K_{0,t}$ and $\mu$ are independent,
$\lambda\otimes \P$-a.s. $K_{0,t}(x)$ is atomic. Using the rotation
invariance, it implies that $K_{0,t}(x)$ is a.s. atomic for all $x$.
\qed

\bibliographystyle{apalike}

\end{document}